\documentclass[12pt,oneside]{amsart}
\usepackage[T1]{fontenc}
\usepackage{geometry}
\usepackage{verbatim}
\usepackage{amssymb}
\usepackage{color}
\makeatletter
 \theoremstyle{plain}
\newtheorem{thm}{Theorem}[section]
  \theoremstyle{plain}
  
  \theoremstyle{remark}
  \newtheorem{rem}[thm]{Remark}
  \theoremstyle{plain}
  
  \theoremstyle{plain}
  \newtheorem{lem}[thm]{Lemma}

\def\bfR#1{{\bf R}^#1}

\def\com#1{ \quad\hbox{#1}\quad}

\smallskip
\def\<{{\langle }}
\def\>{{\rangle }}

\makeatother

\begin{document}

\title{Algebraic zero mean curvature varieties in semi-riemannian manifolds}

\author{ Oscar M. Perdomo }

\date{\today}

\curraddr{Department of Mathematics\\
Central Connecticut State University\\
New Britain, CT 06050\\
}

\email{ perdomoosm@ccsu.edu}

\begin{abstract}

In this paper we provide a family of  algebraic space-like surfaces in the three dimensional anti de Sitter space that shows that this Lorentzian manifold admits algebraic maximal examples of any order.
Then, we classify all the space-like order two algebraic  maximal hypersurfaces in the anti de Sitter $N$-dimensional space. Finally, we provide two families of examples of  Lorentzian order two algebraic zero mean curvature in the de Sitter space.
\end{abstract}

\subjclass[2000]{53C42, 53C50}

\maketitle

\section{Introduction}

Let us denote by $\bfR{{N+1}}_k$ the set $\bfR{{N+1}}$ endowed with the semi-riemannian metric $ds_k^2=-dx_1^2-\dots -dx_k^2+dx_{k+1}^2+\dots+dx_{N+1}^2$ and for  $\epsilon=\pm 1$, let us denote by $K_{k,\epsilon}$ the set of points in $\bfR{{N+1}}_k$ with square norm equal to $\epsilon$, endowed with the metric induced by $ds^2_k$. Notice that  $K_{2,-1}$ is the anti de Sitter $N$ dimensional space, $K_{1,1}$ is the $N$ dimensional de Sitter space, $K_{1,-1}$ is the $N$ dimensional hyperbolic space and $K_{0,1}$ is the Euclidean $N$ dimensional sphere. In this paper we will be dealing with zero mean curvature hypersurfaces in the spaces $K_{k,\epsilon}$. These hypersurfaces are critical points of the area functional and they are called minimal or maximal depending on the nature of the critical point. We will be referring to them as ZMC hypersurfaces.

In 1967, Hsiang studied order $n$ algebraic  minimal hypersurfaces of spheres, \cite{H}, in his paper he showed that such a hypersurface $M\subset S^N$ is defined by the zero level set of an irreducible homogeneous polynomial $f:\bfR{{N+1}}\to {\bf R}$ of degree $n$, if and only if,

\begin{eqnarray}\label{Eq in spheres}
2 \Delta f(x) \, |\nabla f|^2(x)-\<\nabla |\nabla f|^2(x),\nabla f(x)\>=0\com{whenever} f(x)=0
\end{eqnarray}

In the same paper Hsiang proved that the only order two algebraic minimal hypersurfaces in the sphere are the Clifford tori. Recall that Clifford tori have exactly two principal curvatures everywhere and the norm the second fundamental is constant. In 1970, \cite{Lawson}, Lawson proved that the three dimensional sphere admits algebraic surfaces of any order by given explicit examples.

 In this paper we will generalize the results we just mentioned above. First, we will provide a formula similar to (\ref{Eq in spheres}) for order $n$ ZMC algebraic hypersurfaces of $K_{k,\epsilon}$, then, we will classify all order two algebraic ZMC hypersurfaces in the anti de Sitter space. In this classification we prove that the order 2 examples are the hyperbolic cylinders and a family of hypersurfaces with exactly three principal curvatures everywhere and with non constant norm of the second fundamental. We continue the paper by proving that the three dimensional anti de Sitter space admits algebraic surfaces of any order by given explicit examples. All these examples are space-like and complete.

 We will finish the paper by showing two families of algebraic order 2 embedded time-like ZMC hypersurfaces in the $N$ dimensional de Sitter. Hypersurfaces in the first family have three principal curvatures everywhere, $0$ is one of them with multiplicity 1 and corresponds with a time-like direction. Hypersurfaces in the second family have two principal curvatures everywhere, one of them has multiplicity 1 and corresponds with a time-like direction.

We would like to point out that the fact that a space-like hypersurface of a semi riemannian  manifold is closed, does not guarantee the completeness of its induced riemannian metric see e.g, \cite{C-Y}. For this reason, even for algebraic examples which are evidently closed sets, the completeness is a property that must be checked.

\section{Preliminaries}

\medskip

For any non negative integer $s$ less than $N+1$, and any positive integer $j$, we will denote  by $I_j$ the $j\times j$ identity matrix and by $B_s=\{b^s_{ij}\}$ the $(N+1)\times(N+1)$ matrix defined by

\[
B_s=\left(\begin{array}{cc}
-I_s & {\bf 0} \\
{\bf 0} & I_{N+1-s}\end{array}\right),\]

Given $\epsilon = \pm 1$, we will also denote by

$$K_{s,\epsilon}= \{x\in \bfR{{N+1}}: \<B_sx,x\>=\epsilon \} $$

Given an homogeneous polynomial $f:\bfR{{N+1}}\to{\bf R}$, we will define

$$M^{s,\epsilon}_f=\{x\in K_{s,\epsilon}:f(x)=0 \}$$

We will define the metric $g_s:\bfR{{N+1}}\times \bfR{{N+1}}\to {\bf R}$ by $g_s(v,w)=\<B_sv,w\>$, where $\<\, \, ,\, \>$ denotes the Euclidean dot product. Notice that $K_{2,-1}$ with the metric induced by $g_2$ is the Anti de Sitter space, $K_{1,-1}$ with the metric induced by $g_1$ is the Hyperbolic space, and $K_{1,1}$ with the metric induced by $g_1$ is the de Sitter space. In this section we will derive the zero mean curvature equation for the hypersurface $M^{s,\epsilon}_f$. For the sake of completeness we will start this section by proving the following well known property of the Laplacian with respect to a metric $g_s$.

\begin{lem}

\label{Laplacian}

Let  $s$ be any non negative integer less than $N+1$. If for a given smooth function $f:\bfR{{N+1}}\to {\bf R}$ we define

$$\Delta_{g_s}f=-\frac{\partial^2 f}{\partial x_1^2}-\dots-\frac{\partial^2 f}{\partial x_s^2}+\frac{\partial^2 f}{\partial x_{s+1}^2}+\dots+\frac{\partial^2 f}{\partial x_{N+1}^2}$$

then, for any basis $\{v_1,\dots,v_{N+1}\}$ such that $\<B_sv_i,v_j\>=b^s_{ij}$, we have that

$$\Delta_{g_s}f=-\<D^2f(v_1),v_1\>-\dots -\<D^2f(v_s),v_s\>+\<D^2f(v_{s+1}),v_{s+1}\>+\dots+\<D^2f(v_{N+1}),v_{N+1}\>$$

Where $D^2f$ denotes the second derivative $(N+1) \times (N+1)$ matrix of the function $f$.

\end{lem}

\begin{proof}
Let us denote by $e_i=(0,\dots,0,1,0,\dots,0)$, $i=1,\dots ,N+1$ the standard basis of $\bfR{{N+1}}$ and let us define the matrix $A=\{a_{ij}\}$ by the equations $v_i=\sum_{k=1}^{N+1}a_{ik}e_k$. From the condition on the values for $\<B_sv_i,v_j\>$ we get that, $B_s=AB_sA^T$, from this equation we get that the inverse of the matrix $B_sA$ is the matrix $B_sA^T$ and since every matrix commutes with its inverse we get the equation $B_s=A^TB_sA$. We have that

\begin{eqnarray*}
& &-\<D^2f(v_1),v_1\>-\dots -\<D^2f(v_k),v_k\>+\<D^2f(v_{k+1}),v_{k+1}\>+\dots+\<D^2f(v_{N+1}),v_{N+1}\> \\
& &=\sum_{i,j=1}^{N+1} b^s_{ij}\<D^2f(v_i),v_j\> \\
& &= \sum_{i,j,k,l=1}^{N+1} b^s_{ij}a_{ik}a_{jl}\<D^2f(e_k),e_l\> \\
& &= \sum_{k,l=1}^{N+1} b^s_{kl}\<D^2f(e_k),e_l\> \\
& &=\Delta_{g_s}f
\end{eqnarray*}

The third equation follows from the equation $B=A^TBA$.

\end{proof}

\medskip

\begin{rem}
\label{regularity of M}

Recall that for any $x\in  K_{s,\epsilon}\subset \bfR{{N+1}}_s $, a normal vector (with respect to the Euclidean metric) is given by $B_sx$. Therefore, we can describe the tangent space of $ K_{s,\epsilon}$ at $x$ as,

$$T_x K_{s,\epsilon}=\{v\in\bfR{{N+1}}\, :\, \<v,B_sx\>=0  \} =\{v\in\bfR{{N+1}}\, :\, g_s(v,x) = 0 \, \} $$

 Let us assume that $f:\bfR{{N+1}}\to {\bf R}$  is a homogeneous polynomial of degree $k$. Anytime $\nabla f(x)\ne{\bf 0}$ for some $x\in M^{s,\epsilon}_f$, we get that $\< \nabla f(x),x\>=k f(x)=0$. From this equation we get
two important observations. First, we get that in a neighborhood of $x$, $M^{s,\epsilon}_f$ is an embedded $N-1$ submanifold of $\bfR{{N+1}}$, because $\nabla f(x)$ cannot be a multiple of the vector $B_sx$, since the equation $\nabla f(x)=\lambda B_sx$ implies that $0=\< \nabla f(x),x\>=\lambda\<B_sx,x\>=\epsilon \lambda$. The second observation is that the vector $B_s\nabla f(x)\in T_x K_{s,\epsilon}$ and it is perpendicular to $T_xM^{s,\epsilon}_f$ with respect to the metric $g_s$. This follows because $0=\< \nabla f(x),x \>=\< B_s\nabla f(x),B_sx\>=g_s(B_s\nabla f(x),x)$ and for any $v\in T_xM^{s,\epsilon}_f$, $ 0=\< \nabla f(x),v\>=\<B_s \nabla f(x),B_s v\>=g_s(B_s\nabla f(x),v) $. From the two observations above we get that $M^{s,\epsilon}_f$ is an embedded $N-1$ sub manifold near $x$. Moreover, if $g_s(\nabla f(x),\nabla f(x))\ne0$, then, near this point $x$ the metric $g_s$ induces a semi riemannian metric and the vector field $\nu(x)=\frac{1}{\sqrt{|g(B_s\nabla f,B_s\nabla f)|}}B_s\nabla f$ defines a Gauss map in $M^{s,\epsilon}_f\subset K_{s,\epsilon}$.
\end{rem}

\medskip

\begin{rem}
\label{shape operator and mean curvature}

If $M^{s,\epsilon}_f$  is a semi-riemannian hypersurface and $\nu:M^{s,\epsilon}_f\to\bfR{{N+1}}$ is a Gauss map, then, the differential of
$\nu$ at any $x\in M^{s,\epsilon}_f$ defines a linear transformation from $T_xM^{s,\epsilon}_f $ to $T_x M^{s,\epsilon}_f$. This transformation is self-adjoint with respect to the metric $g_s$. We have that, if the metric induced by $g_s$ on $M^{s,\epsilon}_f$ is riemannian, then, the shape operator can be diagonalized, if this metric is not riemannian but only semi riemannian, the shape operator may not be possible to diagonalize. The trace of the shape operator divided by the dimension of $M^{s,\epsilon}_f$ is called the mean curvature and it is denoted by $H(x)$. We have that

$$H(x)=\frac{1}{N-1}\, \sum_{i=1}^{N-1}g_s(v_i,v_i)g_s(d\nu_x(v_i),v_i)$$

where $\{v_1,\dots,v_{N-1}\}$ is any basis for $T_xM^{s,\epsilon}_f$ such that $|g_s(v_i,v_i)|=\delta_{ij}$

\end{rem}

The following lemma gives us a condition in terms of $f$ that guarantees that the hypersurface $M^{s,\epsilon}_f$ has zero mean curvature.

\begin{lem}

\label{The maximal equation}
Let $\epsilon$ be either $1$ or $-1$ and let  $f:\bfR{{N+1}}\to {\bf R}$ be an homogeneous polynomial such that  $w(x)=\<B_s\nabla f(x),\nabla f(x)\>\ne 0$ for every $x$ in an open $U\subset \bfR{{N+1}}$. If

$$\Sigma=\{x\in \bfR{{N+1}}:f(x)=0 \com{and} \<B_sx,x\>=\epsilon \}\cap U$$

is not the empty set, then, $\Sigma$ with the metric induced by $g_s$ is a semi riemannian embedded submanifold and its mean curvature, as a subset of $K_{s,\epsilon}$, vanishes  if and only if

$$ 2 w(x) \Delta_{g_s}f(x) - \<\nabla w(x),B_s \nabla f(x) \>=0\com{for any} x\in \Sigma $$

\end{lem}

\begin{proof}

The first part follows from the Remark (\ref{regularity of M}). From the same remark we also know that the Gauss map is given by $\nu=\frac{B\nabla f}{\sqrt{|w|}}$.

\begin{eqnarray*}
g_s(d\nu_x(v_i),v_i)&=&\<d(\frac{B_s \nabla f}{\sqrt{|w|}})_x(v_i),B_s v_i\>\\
     &=& \<\frac{1}{\sqrt{|w(x)|}} d(B_s \nabla f)_x(v_i),B_s v_i\>+ d(\frac{1}{\sqrt{|w|}} )_x(v_i)\<B_s \nabla f(x),B_sv_i \>\\
&=&\frac{1}{\sqrt{|w(x)|}}\<B_s d(\nabla f)_x(v_i), B_s v_i\>+0\\
&=&\frac{1}{\sqrt{|w(x)|}} \<D^2f(x)v_i,v_i\>\\
\end{eqnarray*}

On the other hand, since $\<\nabla f(x),x\>=kf(x)$, we get that for every $x\in \Sigma$, $\<D^2f(x)x,x\>=0$. We also have that the derivative of the function $w$ in the direction $B_s\nabla f(x)$ can be computed as follows,

\begin{eqnarray*}
dw_x(B_s\nabla f(x))&=& B_s\nabla f(x)\<\nabla f,B_s\nabla f\> \\
&=&\<d(\nabla f)_x(B_s\nabla f(x)),B_s\nabla f(x)\>+\<\nabla f(x),d( B_s  \nabla f)_x(B_s\nabla f(x))\>\\
&=&\<D^2f(x) B_s\nabla f(x),B_s\nabla f(x)\>+\<\nabla f(x), B_s  d( \nabla  f)_x(B_s\nabla f(x))\>\\
&=&\<D^2f(x) B_s\nabla f(x),B_s\nabla f(x)\>+\<B_s \nabla f(x),   D^2f(x)(B_s\nabla f(x))\>\\
&=&2 \<D^2f(x) B_s\nabla f(x),B_s\nabla f(x)\>\\
\end{eqnarray*}

Now, if $\{v_1,\dots,v_{N-1}\}$ is a basis of $T_x\Sigma$ with $|\<B_sv_i,v_j\>|=\delta_{ij}$ and $\hbox{sig}(a)$ denotes the sign of the number $a$, we get that

\begin{eqnarray*}
\Delta_{g_s}f &=&\epsilon \, \<D^2f(x)\, x,x\>+\hbox{sig}(w(x))\<D^2f(x)\nu(x),\nu(x)\>+\sum_{i=1}^{N-1}\hbox{sig}(\<B_sv_i,v_i\>)\<D^2f(x)v_i,v_i\>  \\
&=&  \hbox{sig}(w(x))\frac{1}{|w(x)|^2}\<D^2f(x)B_s\nabla f(x),B_s\nabla f(x)\>+\sum_{i=1}^{N-1}\hbox{sig}(\<B_sv_i,v_i\>)\<D^2f(x)v_i,v_i\>  \\
&=&  \frac{1}{2w(x)} dw_x(B_s\nabla f(x)) +\sum_{i=1}^{N-1}\hbox{sig}(\<B_sv_i,v_i\>)\<D^2f(x)v_i,v_i\>  \\
&=&  \frac{1}{2w(x)}\<\nabla w(x),B_s\nabla f(x)\> +(N-1)\, H(x) \\
\end{eqnarray*}

The last equation implies the lemma.

\end{proof}

\medskip
{\bf Definition:} {\sl
We say that a set $M\subset K_{s,\epsilon}$ is algebraic, if there exists an irreducible homogeneous polynomial $f:\bfR{{N+1}}\to {\bf R}$ such that
$M\subset f^{-1}(0)$ and

$$2 w(x) \Delta_{g_s}f(x) - \<\nabla w(x),B_s \nabla f(x) \>=0\com{for every} x\in M$$

where $w=\<\nabla f, B_s\nabla f\>$. We will say that an algebraic hypersurface $M\subset K_{s,\epsilon}$ is regular, if $w(x)\ne0$ for every $x\in M$.
}

\medskip

\begin{rem}
\label{dual hypersurfaces}
If a homogeneous polynomial $f:\bfR{{N+1}}\to {\bf R}$ satisfies that
\begin{eqnarray}
\label{algebraic equation} 2 w(x) \Delta_{g_s}f(x) - \<\nabla w(x),B_s \nabla f(x) \>=h(x)f(x)
\end{eqnarray}

then, the function $f$ defines zero mean curvature hypersurfaces in $K_{s,\epsilon}$ for either $\epsilon=1$ or $\epsilon=-1$. In particular,
homogeneous polynomials satisfying equation (\ref{algebraic equation}) for $s=1$ define algebraic minimal hypersurfaces in the Hyperbolic space and
algebraic zero mean curvature hypersurfaces in the de Sitter space.
\end{rem}

It is not known if in general, a polynomial that defines an algebraic hypersurface, satisfies the equation ( \ref{algebraic equation}) for some homogeneous polynomial $h$.


{\bf Conjecture: } {\sl If $M\subset K_{s,\epsilon}$ is an algebraic hypersurface defined by a homogeneous polynomial $f$, then, there exists a homogeneous
polynomial $h$ such that
\begin{eqnarray}
2 w(x) \Delta_{g_s}f(x) - \<\nabla w(x),B_s \nabla f(x) \>=h(x)f(x)
\end{eqnarray}
}
\medskip

\section{Algebraic surfaces in the anti de Sitter space}

In this section we will show explicit examples  of algebraic surfaces of any order in the anti de Sitter space.

\medskip

\begin{thm}
 Let $k$ and $n$ be any pair of positive relative prime integers such that $n$ is odd.

 \medskip

{\bf a.} If $k<n$,  the map $\phi:\bfR{2}\to K_{2,-1}$ given by

$$\phi(s,t)=(\cosh(s) \cosh(n t), \sinh(s) \sinh(k t),
\cosh(s) \sinh(n t), -\cosh(k t) \sinh(s)) $$

defines an order $k+n$ algebraic ZMC space-like complete immersed surface in the anti de Sitter space $K_{2,-1}$.
 \medskip

{\bf b.} If $k>n$, the map $\rho:\bfR{2}\to K_{2,-1}$ given by

$$\rho(s,t)=(\cosh(n t) \sinh(s), \cosh(s) \sinh(k t),
\sinh(s) \sinh(n t), -\cosh(s) \cosh(k t)) $$

defines an order $k+n$ algebraic ZMC immersed surface in  $K_{2,1}$. Moreover, $g_2$ induces a negative  definite semi-riemannian  metric.

{\bf c.} The polynomial

$$f=2 \, ((x_1 - x_3)^k (x_2 - x_4)^n + (x_1 + x_3)^k (x_2 + x_4)^n) $$

defines algebraic minimal surfaces in the spaces $K_{2,1}$  and $K_{2,-1}$. Moreover, $f(\phi(s,t))=0$ and  $f(\rho(s,t))=0$ for every $s$ and $t$.

\end{thm}

\medskip

\begin{proof}

Let us start proving {\bf c.} A direct verification shows that if we define
$$ w=\<B_2 \nabla f,\nabla f\>\com{ and} g=2 w \Delta_{g_2}f-\<B_2\nabla f,\nabla w\>$$

then $g=hf$ where

\begin{eqnarray*}
h&=&-32 k^2 (x_1^2 - x_3^2)^{k-2} (2 x_1^2 - 2 x_3^2 + k^2 (x_2^2 - x_4^2) + k (-x_2^2 + x_4^2))\com{when} n=1\\
h&=& -32 n^2 (x_2^2 - x_4^2)^{n-2} (n^2 (x_1^2 - x_3^2) + n (-x_1^2 + x_3^2) + 2 (x_2^2 - x_4^2))\com{when} k=1
\end{eqnarray*}

and

\begin{eqnarray*}
h&=&-32 \, (x_1^2 - x_3^2)^{k-2} (\, x_2^2 - x_4^2)^{
  n-2} ( n^4 (x_1^2 - x_3^2)^2-n^3 (x_1^2 - x_3^2)^2  +\\
   & & 2 k^2 n^2 (x_1^2 - x_3^2) (x_2^2 - x_4^2) +
   (k-1) k^3 (x_2^2 -x_4^2)^2\, ) \com{for $n\ge 2$ and $k\ge2$}
\end{eqnarray*}

Therefore, the function $f$ induces algebraic examples in the spaces $K_{2,1}$  and $K_{2,-1}$. A direct computation shows that,

 $$f(\phi(s,t))=2 \cosh^k(s)\sinh^n(s)\, (1+(-1)^n)$$

 and

 $$f(\rho(s,t))=2 \cosh^n(s)\sinh^k(s)\, (1+(-1)^n)$$

 since $n$ is odd, then, $f(\phi(s,t))=0$ and  $f(\rho(s,t))=0$ for every $s$ and $t$. A direct verification shows that $\<B_2\phi(s,t),\phi(s,t)\>=-1$ and $\<B_2\rho(s,t),\rho(s,t)\>=1$, therefore, $\phi(s,t)\in K_{2,-1}$ and
 $\rho(s,t)\in K_{2,1}$.  The fact that the map $\phi$ is an immersion follows from the following equations,

\begin{eqnarray*}
\<\frac{\partial \phi}{\partial s},B_2\frac{\partial \phi}{\partial s}\>&=&  1\com{and} \\
\<\frac{\partial \phi}{\partial s},B_2\frac{\partial \phi}{\partial t}\> &=&0\\
\<\frac{\partial \phi}{\partial t},B_2\frac{\partial \phi}{\partial t}\> &=&\frac{1}{2} (\, k^2 + n^2 + (n^2-k^2 ) \cosh(2 s)\, )
\end{eqnarray*}

Therefore, since $n>k$, $\phi$ is a space-like complete immersion. Notice that we already proved the ZMC condition because $f(\phi(s,t))=0$.
The fact that the map $\rho$ is an immersion follows from the following equations,

\begin{eqnarray*}
\<\frac{\partial \rho}{\partial s},B_2\frac{\partial \rho}{\partial s}\>&=&  -1\com{and} \\
\<\frac{\partial \rho}{\partial s},B_2\frac{\partial \rho}{\partial t}\> &=&0\\
\<\frac{\partial \rho}{\partial t},B_2\frac{\partial \rho}{\partial t}\> &=&-\frac{1}{2} (\, k^2 + n^2 + (k^2-n^2 ) \cosh(2 s)\, )
\end{eqnarray*}

Therefore, since $k>n$, $\rho$ is an immersion and the metric induced by $g_2$ is negative definite. Notice that we already proved the  ZMC condition because $f(\rho(s,t))=0$.

\end{proof}

\begin{rem}

The reason we ask for $n$ and $k$ to be relative primes is that if they are not, then the polynomial $f$ is not irreducible.

The existence of algebraic surfaces by pairs, in $H_{2,1}$ and $H_{2,-1}$, is a consequence of remark (\ref{dual hypersurfaces}).

\end{rem}

\section{\label{examples anti de Sitter} Maximal hypersurfaces in  anti de Sitter spaces}

\begin{thm}
For any pair of positive integers $m$ and $n$ and any non negative integer $k$ the polynomial $f:\bfR{{2+m+n+k}}\to{\bf R}$ defined by

$$ f(x_1,x_2,y,z,u) =2\, x_1 x_2 + \frac{m - n}{\sqrt{m n}}\,  x_2^2 + \sqrt{\frac{n}{m}}\,  |y|^2 - \sqrt{\frac{m}{n}}\, |z|^2$$

where $y\in \bfR{m}$, $z\in \bfR{n}$ and $u\in \bfR{k}$, provides complete ZMC space-like regular algebraic varieties in the anti de Sitter space. If we denote by  $\Sigma = f^{-1}(0)\cap K_{2,-1}$ the hypersurface induced by $f$, we have the following properties: if $k=0$, then, $\Sigma$ is a hyperbolic cylinder with principal curvatures $-\sqrt{\frac{n}{m}}$ and $\sqrt{\frac{m}{n}}$ with multiplicity $m$ and $n$ respectively. If $k>0$, then, the principal curvatures at a point $(x_1,x_2,y,z,u)\in\Sigma$ are $\, 0$, $-\sqrt{\frac{n}{m(1+|u|^2)}}$ and $\sqrt{\frac{m}{n(1+|u|^2)}}$ with multiplicities $k$, $m$ and $n$ respectively.

\end{thm}

\begin{proof}

Let us define $w=\<\nabla f, B_2\nabla f\>$.
A direct computation shows that

$$w=-4 x_2^2 - 4 \big{(}\, x_1 + \frac{(m - n) x_2}{\sqrt{m n}}\, \big{)}^2 + 4 \frac{n}{m} |y|^2 + 4 \frac{m}{n} |z|^2\com{and}\Delta_{g_2}f=\frac{2(n-m)}{\sqrt{nm}}$$

Using this expression for $w$ we get that $\<\nabla w,B_2 \nabla f\>$ reduces to

$$16\, \big{(}\, \frac{(m-n)(m^2+n^2)}{(mn)^{\frac{3}{2}}} x_2^2 +\frac{2(m^2-mn+n^2)}{mn} x_1x_2+\frac{m-n}{\sqrt{mn}} x_1^2+\sqrt{\frac{n^3}{m^3}}\, |y|^2-\sqrt{\frac{m^3}{n^3}}\, |z|^2\,\big{)}$$

Finally, if we define $g=2 w(x) \Delta_{g_2}f(x) - \<\nabla w(x),B_2 \nabla f(x) \>$, a direct computation shows that $g=-16 f$. Therefore $f$ defines an algebraic ZMC hypersurface. We will prove that $\Sigma$ is an embedded hypersurface by showing that $w$ never vanishes for points $p=(x_1,x_2,x_3,y,z)$ in $\Sigma$. Solving the equations $f=0$ and $\<B_2p,p\>=-1$ for $|y|^2$ and $|z|^2$ we get that

\begin{eqnarray*}
|\, y|^2 &=& \frac{(\sqrt{m}\, x_1-\sqrt{n}\, x_2)^2-m(1 + |u|^2\, )}{
 m + n}  \\
 |\, z|^2 &=& \frac{(\sqrt{m}\, x_2+\sqrt{n}\, x_1)^2-n(1 + |u|^2\, )}{
 m + n}
\end{eqnarray*}

Substituting these values for $|\, y|^2$ and $|z|^2$ in our previous expression for $w$ we get that

$$w=-4\, (1+|u|^2)$$

Therefore $\Sigma $ is an embedded hypersurface. As pointed out in Remark (\ref{regularity of M}), the Gauss map of $\Sigma$ is given by $\nu=\frac{1}{\sqrt{|w|}}B_2\nabla f$. For points in $\Sigma$, $\nu$ reduces to

$$\nu=\frac{1}{\sqrt{ 1+|u|^2}}\, (\, -x_2,   \frac{n - m}{\sqrt{m n}}\, x_2-x_1 ,  \sqrt{\frac{n}{m}} \, y, -\sqrt{\frac{m}{n}} \, z, {\bf 0} \,) $$

A direct computation shows that if   $p=(x_1,x_2,y,z,u)\in\Sigma$, then the subspace $\Gamma$ defined by,

$$ \{ (0,0,v,{\bf 0},{\bf 0})+\lambda\, (\sqrt{m},-\sqrt{n},(\sqrt{m}\, x_1-\sqrt{n}\, x_2)\, \frac{y}{|y|^2},{\bf 0},{\bf 0}\, ) \in\bfR{{2+m+n+k}}\, :\, v\in\bfR{m}\, , \<v,y\>=0\, ,\lambda\in {\bf R}\,\} $$

in the case that $|\, y|\ne0$, or by

$$ \{ (0,0,v,{\bf 0},{\bf 0})\in\bfR{{2+m+n+k}}\com{where}v\in\bfR{m}\, \} \com{when $|\, y|=0$}$$

is contained in $T_p\Sigma$ and $d\nu_p(v)=\sqrt{\frac{n}{m(1+|u|^2)}} \, v$ for any $v\in \Gamma$. Notice that the dimension of $\Gamma$ is $m$. Also, notice that when $k=0$ the expression $-\sqrt{\frac{n}{m(1+|u|^2)}}$ reduces to $-\sqrt{\frac{n}{m}}$.

Likewise, for the same point $p$, a direct computation shows that the subspace $\Omega$ defined by,

$$ \{ (0,0,{\bf 0},v,{\bf 0})+\lambda\, (\sqrt{n},\sqrt{m},0,(\sqrt{m}\, x_2+\sqrt{n}\, x_1)\, \frac{z}{|z|^2},{\bf 0}\, ) \in\bfR{{2+m+n+k}}\, :\, v\in\bfR{n}\, , \<v,z\>=0\, ,\lambda\in {\bf R}\,\} $$

in the case that $|\, z|\ne0$ or by

$$ \{ (0,0,{\bf 0},v,{\bf 0})\in\bfR{{2+m+n+k}}\com{where}v\in\bfR{n}\, \} \com{when $|\, z|=0$}$$

is contained in $T_p\Sigma$ and $d\nu_p(v)=-\sqrt{\frac{m}{n(1+|u|^2)}} \, v$ for any $v\in\Omega$. Notice that the dimension of $\Omega$ is $n$. Also, notice that when $k=0$ the expression $\sqrt{\frac{m}{n(1+|u|^2)}}$ reduces to $\sqrt{\frac{m}{n}}$.

Also, a direct computation shows that,  when $k>0$, the vector space $\Pi$ given by

$$ \{ (0,0,{\bf 0},{\bf 0},v)+\frac{\lambda |u|^2}{1+|u|^2}\, ( x_1, x_2, y, z,\frac{1+|u|^2}{|u|^2}\, u\, ) \in\bfR{{2+m+n+k}}\, :\, v\in\bfR{k}\, , \<v,u\>=0\, ,\lambda\in {\bf R}\,\} $$
in the case that $|\, u|\ne0$ or by

$$ \{ (0,0,{\bf 0},{\bf 0},v)\in\bfR{{2+m+n}}\com{where}v\in\bfR{k}\, \} \com{when $|\, u|=0$}$$

is contained in $T_p\Sigma$ and $d\nu_p(v)=0$ for any $v\in\Pi$. Therefore, since the dimension of $\Sigma$ is $\hbox{dim}(\Gamma)+\hbox{dim}(\Omega)+\hbox{dim}(\Pi)$, we conclude that  $0$ is a principal curvature with multiplicity $k$,  $-\sqrt{\frac{n}{m(1+|u|^2)}} $ is a principal curvature with multiplicity $m$, and $\sqrt{\frac{m}{n(1+|u|^2)}}$ is a principal curvature with multiplicity $n$. The completeness of these examples follows because these examples are either isometry to hyperbolic cylinders or they are isometry to $M\times  \bfR{k}$ endowed with the metric $ds^2(p,v)=(1+|v|^2) ds_1^2(p)+ds_2^2(v)$ where $(M,ds_1^2)$ is a hyperbolic cylinder, and $(\bfR{k},ds_2^2)$ is the Euclidean space with the metric $ds_2^2$ defined by the $k\times k$ symmetric matrix

$$D(v)=I-\frac{1}{1+|v|^2}v\, v^T\com{where } v^T=(v_1,\dots, v_k)$$

 $(\bfR{k},ds_2^2)$ is complete because the eigenvalues of the matrix $D(v)$ are $\frac{1}{1+|v|^2}$  and $1$. Notice that the eigenvalue of $D(v)$ are either vectors parallel to $v$ or vectors perpendicular to the vector $v$. The isometry mentioned above is given by the map $h:M\times  \bfR{k}\to \Sigma$ defined as

$$h(p,v)=(\sqrt{1+|v|^2}\, p\, ,\, v\,)$$

\end{proof}

\medskip

\section{ Zero mean curvature hypersurfaces in the de Sitter space}

\begin{thm}
For any pair of positive integers $m$ and $n$ the polynomial $f:\bfR{{3+m+n}}\to{\bf R}$ defined by

$$ f(x_1,x_2,x_3,y,z) = 2 x_3 x_2 +\frac{n-m}{\sqrt{mn}} x_2^2 + \sqrt{\frac{n}{m}} |y|^2 -\sqrt{\frac{m}{n}} |z|^2$$

where $y\in \bfR{m}$ and $z\in \bfR{n}$ provides zero mean curvature regular algebraic  varieties in the de Sitter space. Moreover, for any $p=(x_1,x_2,x_3,y,z)\in \Sigma=f^{-1}(0)\cap K_{1,1}$, the principal curvatures are $0$, $-\sqrt{\frac{n}{m(1+x_1^2)}}$ and $\sqrt{\frac{m}{n(1+x_1^2)}}$ with multiplicity $1$, $m$ and $n$ respectively. Additionally,  $\Sigma$ with the metric induced by $g_1$ is a Lorentzian semi riemannian and the principal direction associated with the principal curvature $0$ is a time-like direction.

\end{thm}

\begin{proof}

Let us define $w=\<\nabla f, B_1\nabla f\>$.
A direct computation shows that

$$w=4\big{(}\,  x_2^2 +   \frac{((n - m) x_2+\sqrt{mn}\, x_3)^2}{m n}\, +  \frac{n}{m} |y|^2 +  \frac{m}{n} |z|^2\big{)}\com{and}\Delta_{g_1}f=\frac{2(n-m)}{\sqrt{nm}}$$

Using this expression for $w$ we get that $\<\nabla w,B_1 \nabla f\>$ reduces to

$$16\, \big{(}\, \frac{(n-m)(m^2+n^2)}{(mn)^{\frac{3}{2}}} x_2^2 +\frac{2(m^2-mn+n^2)}{mn} x_3x_2+\frac{n-m}{\sqrt{mn}} x_3^2+\sqrt{\frac{n^3}{m^3}}\, |y|^2-\sqrt{\frac{m^3}{n^3}}\, |z|^2\,\big{)}$$

Finally, if we define $g=2 w(x) \Delta_{g_1}f(x) - \<\nabla w(x),B_1 \nabla f(x) \>$, a direct computation shows that $g=-16 f$. Therefore, $f$ defines an algebraic ZMC hypersurface. We will prove that $\Sigma$ is an embedded hypersurface by showing that $w$ never vanishes for points $p=(x_1,x_2,x_3,y,z)$ in $\Sigma$, solving the equations $f=0$ and $\<B_1p,p\>=1$ for $|y|^2$ and $|z|^2$ we get that

\begin{eqnarray*}
|\, y|^2 &=& \frac{m(1 + x_1^2\, )-(\sqrt{m}\, x_3+\sqrt{n}\, x_2)^2}{
 m + n}  \\
 |\, z|^2 &=& \frac{n(1 + x_1^2\, )-(\sqrt{m}\, x_2-\sqrt{n}\, x_3)^2}{
 m + n}
\end{eqnarray*}

Substituting these values for $|\, y|^2$ and $|z|^2$ in our previous expression for $w$ we get that

$$w=4\, (1+x_1^2)$$

Therefore, $\Sigma $ is an embedded Lorentzian  hypersurface. As pointed out in Remark (\ref{regularity of M}), the Gauss map of $\Sigma$ is given by $\nu=\frac{1}{\sqrt{|w|}}B_1\nabla f$. For points in $\Sigma$, $\nu$ reduces to

$$\nu=\frac{1}{\sqrt{ 1+x_1^2}}\, (\,0,   \frac{n - m}{\sqrt{m n}}\, x_2+x_3 ,x_2 , \sqrt{\frac{n}{m}} \, y, -\sqrt{\frac{m}{n}} \, z\, ) $$

A direct computation shows that if   $p=(x_1,x_2,x_3,y,z)\in\Sigma$, then the subspace $ \Gamma$ defined by,

$$\{ (0,0,0,v,{\bf 0},{\bf 0})+\lambda\, (0,\sqrt{n},\sqrt{m},-(\sqrt{m}\, x_3+\sqrt{n}\, x_2)\, \frac{y}{|y|^2},{\bf 0}\, ) \in\bfR{{3+m+n}}\, :\, v\in\bfR{m}\, , \<v,y\>=0\, ,\lambda\in {\bf R}\,\} $$

in the case that $|\, y|\ne0$ or by

$$ \{ (0,0,0,v,{\bf 0},{\bf 0})\in\bfR{{3+m+n}}\com{where}v\in\bfR{m}\, \} \com{when $|\, y|=0$}$$

is contained in $T_p\Sigma$ and $d\nu_p(v)=\sqrt{\frac{n}{m(1+x_1^2)}} \, v$ for any $v\in  \Gamma$. Notice that the dimension of $\Gamma$ is $m$.
Likewise, for the same point $p$, a direct computation shows that the subspace $\Omega$ defined by,

$$ \{ (0,0,0,{\bf 0},v)+\lambda\, (0,\sqrt{m},-\sqrt{n},0,(\sqrt{n}\, x_3-\sqrt{m}\, x_2)\, \frac{z}{|z|^2}\, ) \in\bfR{{3+m+n}}\, :\, v\in\bfR{n}\, , \<v,z\>=0\, ,\lambda\in {\bf R}\,\} $$

in the case that $|\, z|\ne0$ or by

$$ \{ (0,0,0,{\bf 0},v\, )\in\bfR{{3+m+n}}\com{where}v\in\bfR{n}\, \} \com{when $|\, z|=0$}$$

is contained in $T_p\Sigma$ and $d\nu_p(v)=-\sqrt{\frac{m}{n(1+x_1^2)}} \, v$ for any $v\in\Omega$. Notice that the dimension of $\Omega$ is $n$. Also, a direct computation shows that,   the vector

$$ v=  (1+x_1^2\, ,x_1 x_2,x_1 x_3,x_1 y,x_1 z )  $$

is in $T_p\Sigma$ and $d\nu_p(v)=0$. Therefore, since the dimension of $\Sigma$ is $1+\hbox{dim}(\Gamma)+\hbox{dim}(\Omega)$, we conclude that  $0$ is a principal curvature with multiplicity $1$,  $-\sqrt{\frac{n}{m(1+x_1^2)}} $ is a principal curvature with multiplicity $m$ and $\sqrt{\frac{m}{n(1+x_1^2)}}$ is a principal curvature with multiplicity $n$. Finally, the statement  about  the principal direction associated with the zero eigenvalue follows from the equation $\<v,B_1v\>=-1-x_1^2$ for points in $K_{1,1}$.

\end{proof}

\medskip

\begin{thm}
For any positive integer $m$,  the polynomial $f:\bfR{{3+m}}\to{\bf R}$ defined by

$$ f(x_1,x_2,x_3,y) = \sqrt{m}\,  x_1^2 + 2 x_2 x_3-
 \frac{m-1}{\sqrt{m}}\, x_3^2 + \frac{1}{\sqrt{m}}\, |y|^2 $$

where $y\in \bfR{m}$ provides zero mean curvature regular algebraic  varieties in the de Sitter space. Moreover, for any $p=(x_1,x_2,x_3,y,z)\in \Sigma=f^{-1}(0)\cap K_{1,1}$, the principal curvatures are $\sqrt{m}$ and $-\frac{1}{\sqrt{m}}$ and with multiplicity $1$ and $m$ respectively. Additionally,  $\Sigma$ with the metric induced by $g_1$ is a Lorentzian semi riemannian and the principal direction associated with the principal curvature $\frac{1}{\sqrt{m}}$ is a time-like direction.

\end{thm}

\begin{proof}

Let us define $w=\<\nabla f, B_1\nabla f\>$.
A direct computation shows that
$$w=4( (x_2 - \frac{m-1}{\sqrt{m}}\, x_3)^2 - m x_1^2 +  x_3^2  + \frac{1}{m}\, |y|^2)\com{and}\Delta_{g_1}f=\frac{2(1-m)}{\sqrt{m}}$$

Using this expression for $w$ we get that $\<\nabla w,B_1 \nabla f\>$ reduces to

$$16\, \big{(}\, m^{\frac{3}{2}}\, x_1^2 + \frac{1-m}{\sqrt{m}}\, x_2^2+
\frac{2(m^2-m+1)}{m}\, x_2 x_3 +\frac{(1-m)(1+m^2)}{m^{\frac{3}{2}}}\, x_3^2 +\frac{1}{m^{\frac{3}{2}}}\, |y|^2\big{)}$$

Finally, if we define $g=2 w(x) \Delta_{g_1}f(x) - \<\nabla w(x),B_1 \nabla f(x) \>$, a direct computation shows that $g=-16 f$. Therefore, $f$ defines an algebraic ZMC hypersurface. We will prove that $\Sigma$ is an embedded hypersurface by showing that $w$ never vanishes for points $p=(x_1,x_2,x_3,y)$ in $\Sigma$. Solving the equations $f=0$ and $\<B_1p,p\>=1$ for $|y|^2$ and $x_1^2$ we get that

\begin{eqnarray*}
|\, y|^2 &=& \frac{m-(\sqrt{m}\, x_2+ x_3)^2}{
 m + 1}  \\
 x_1^2 &=& \frac{(\sqrt{m}\, x_3- x_2)^2-1}{
 m + 1}
\end{eqnarray*}

Substituting these values for $|\, y|^2$ and $x_1^2$ in our previous expression for $w$ we get that

$$w=4$$

Therefore $\Sigma $ is an embedded Lorentzian  hypersurface. As pointed out in (\ref{regularity of M}), the Gauss map of $\Sigma$ is given by $\nu=\frac{1}{\sqrt{|w|}}B_1\nabla f$. For points in $\Sigma$, $\nu$ reduces to

$$\nu=(-\sqrt{m}\, x_1,x_3, x_2 -  \frac{m-1}{\sqrt{m}}\, x_3 , \sqrt{\frac{1}{m}} \, y) $$

A direct computation shows that if   $p=(x_1,x_2,x_3,y)\in\Sigma$, then the subspace defined $\Gamma$ by,

$$ \{ (0,0,0,v)+\lambda\, (0,\sqrt{m},1,-(\sqrt{m}\, x_2+ x_3)\, \frac{y}{|y|^2}\, ) \in\bfR{{3+m}}\, :\, v\in\bfR{m}\, , \<v,y\>=0\, ,\lambda\in {\bf R}\,\} $$

in the case that $|\, y|\ne0$ or by

$$\{ (0,0,0,v\,)\in\bfR{{3+m+n}}\com{where}v\in\bfR{m}\, \} \com{when $|\, y|=0$}$$

is contained in $T_p\Sigma$ and $d\nu_p(v)=\sqrt{\frac{1}{m}} \, v$ for any $v\in \Sigma$. Notice that the dimension of $\Gamma$ is $m$.
Also, a direct computation shows that the vector

$$ v=  (x_2-\sqrt{m}\, x_3 ,x_1,-\sqrt{m}\, x_1,{\bf 0}\, )$$

satisfies $\<v,B_1 v\>=-1$ for points in $\Sigma$ and, moreover, we have that $v$ is in $ T_p\Sigma$ and $d\nu_p(v)=-\sqrt{m}\, v$. Therefore, since the dimension of $\Sigma$ is $1+\hbox{dim}(\Gamma)$, we conclude that  $\sqrt{m}$ is a principal curvature with multiplicity $1$ and  $-\frac{1}{\sqrt{m}} $ is a principal curvature with multiplicity $m$.

\end{proof}

\section{Classification of order two algebraic space-like maximal varieties in the anti de Sitter space.}

In this section we will prove  that the examples that we studied in section (\ref{examples anti de Sitter})  constitute all algebraic maximal space-like hypersurfaces in the space in $N$-dimensional anti de Sitter space.

\begin{rem}
In the proof of the next two theorems everytime we mention a change of coordinates, we will be referring to a change of coordinates that corresponds to an isometry of the ambient space.

\end{rem}

\begin{thm}
If $f:\bfR{{N+1}}\to {\bf R}$ is a homogeneous polynomial of degree 2 and $M_f=\{x\in \bfR{{N+1}}:f(x)=0 \com{and} \<B_2x,x\>=-1 \}$ is a regular algebraic ZMC space-like variety of the anti de Sitter $N$-dimensional space, then, up to an isometry on the anti de Sitter space, there exist positive integers $m$ and $n$ and a non negative integer $k$ such that

$$  f(x)=2 x_1 x_2 +  \frac{m-n}{\sqrt{m n}} \, x_2^2  +\sqrt{\frac{n}{m}} \, |y|^2-\sqrt{ \frac{m}{n}} \, |z|^2  $$

where $m+n+k=N-1$ and $y=(x_3,\dots,x_{m+2})$ and $z=(x_{m+3},\dots,x_{m+n+2})$

\end{thm}

\begin{proof}
 Let us assume that $f=\<Ax,x\>$ where $A$ is a $N+1$ by $N+1$ symmetric matrix. Let us denote by

$$\{{\bf e_1}=(1,0,\dots,0),\dots,{\bf e_{N+1}}=(0,0,\dots,1)\}$$

the standard basis of $\bfR{{N+1}}$. Since $M_f$ is space-like, without loss of generality, we will assume that the point ${\bf e_1}$  is in $M_f$, $\nabla f({\bf e_1})=2\, {\bf e_2}$, that the principal directions at ${\bf e_1}$ are given by the vectors ${\bf e_3},\dots,{\bf e_{N+1}}$, and that the principal curvatures at ${\bf e_1}$ are $\kappa_1,\dots,\kappa_{N-1}$ respectively. We will additionally assume that $\kappa_1,\dots,\kappa_m$ are positive $\kappa_{m+1},\dots,\kappa_{m+n}$ are negative and $\kappa_{m+n+1}\dots\kappa_{m+n+k}=0$.  Under these assumptions the matrix $A$ reduces to

\[
A=\left(\begin{array}{cccccc}
0 & 1       & 0         & 0   &   \dots             &0           \\
1 & b       & a_1       & a_2 &   \dots             &a_{N-1}     \\
0 & a_1     & \kappa_1  & 0   &   \dots             & 0          \\
0 & a_2     & 0         & \kappa_2 &  0             &0           \\
0 & \vdots  & 0         & 0        &\ddots          &0           \\
0&a_{N-1}   &0          &0         &\dots           & \kappa_{N-1}\end{array}\right).\]

and the function $f$ reduces to

$$f(x)=2 x_1 x_2 + b x_2^2 + 2 a_1 x_2 x_3  + \dots + 2 a_{N-1} x_2 x_{N+1} + \kappa_1 x_3^2+\dots+ k_{N-1} x_{N+1}^2 $$

The ZMC condition on $M_f$ implies that
\begin{eqnarray}\label{principal curvatures}
\kappa_1+\kappa_2+\dots+\kappa_{N-1}=0
\end{eqnarray}

Let us define,

$$w=\<\nabla f, B_2\nabla f\>\com{and} g= 2 w \Delta_{g_2}f - \<\nabla w,B_2 \nabla f \>$$

By Lemma (\ref{The maximal equation}) we have that $g$ vanishes whenever $f$ vanishes. Therefore, for any point in the $M_f$, the gradient of $g$, $\nabla g$ must be parallel to the gradient of $f$, $\nabla f$. A direct computation shows that

$$\nabla g({\bf e_1})=({0,
 32 (-1 + a_1^2 +\dots + a_{N-1}^2 )\, , 32 a_1 \kappa_1,\dots, 32 a_{N-1} \kappa_{N-1}})
$$

In the computation above we have used equation (\ref{principal curvatures}).

Since $\nabla f({\bf e_1})=2 {\bf e_2}$, we get that

$$a_1 \kappa_1=0,\quad a_2 \kappa_2=0,\quad \dots\quad a_{N-1} \kappa_{N-1}=0 $$

The equation above implies that $a_1=0,\dots,a_{m+n}=0$. The function $f$ then reduces to

$$f(x)=2 x_1 x_2 + b x_2^2  + \kappa_1 x_3^2+\dots+ k_{m+n} x_{m+n+2}^2 + 2\, x_2 (a_{m+n+1}  x_{m+n+3}  + \dots + 2 a_{N-1} x_{N+1}) $$

Let us denote by $a=(a_{m+n+1},\dots, a_{m+n+k})$.

Recall that not all the $\kappa_i$'s can be zero, because this would imply that $f$ is  reducible. Therefore we have that $m>0$ and $n>0$. A direct computation shows that for every $j=1,\dots,n$, the point

$$p_1=\sqrt{2}\, {\bf e_1}+\sqrt{\frac{-\kappa_{m+j}}{\kappa_1-\kappa_{m+j}}}\, {\bf e_3}+\sqrt{\frac{\kappa_{1}}{\kappa_1-\kappa_{m+j}}}\,  {\bf e_{m+j+2}}$$

is a point in $M_f$, therefore $g(p_1)$ must be zero. A direct computation shows that

$$g(p_1)=16 \kappa_1 \kappa_{m+j}\, (b + \kappa_1 + \kappa_{m+j})$$

Therefore, $b=-(\kappa_1 + \kappa_{m+j})$. Since this is true for all $j$,  we get that $\kappa_{m+j}=\kappa_{m+1}$ for all $j=1,\dots,n$. In the same way it follows that $\kappa_{i}=\kappa_{1}$ for all $i=1,\dots,m$. Using the equation (\ref{principal curvatures}) we get that

$$\kappa_{m+1}=-\frac{m}{n}\kappa_1 \com{and} b=-\kappa_1-\kappa_{m+1}=\kappa_1\frac{m-n}{n} $$

Denoting by $y=(x_3,\dots,x_{m+2})$ and $z=(x_{m+3},\dots,x_{m+n+2})$, we can write the function $f$ as

$$f(x)=2 x_1 x_2 + \kappa_1\, (\frac{m-n}{n} \, x_2^2  + |y|^2- \frac{m}{n} \, |z|^2 ) + 2 x_2 (a_{m+n+1}  x_{m+n+3}  + \dots + 2 a_{N-1} x_{N+1}) $$

A direct computation shows that the point $p_2=\sqrt{\frac{m(mn + 1+n^2)}{n}}\, {\bf e_2}+\sqrt{mn}\,  {\bf e_{3}}+\sqrt{\frac{-n + n m^2 + m}{n}}\, {\bf e_{m+3}}$ is a point in $M_f$, therefore $g(p_2)$ must be zero. We can check that,

$$g(p_2)=\frac{16 \kappa_1 m (m - n) (\kappa_1^2 m + n( |\, a|^2 - 1) ) (1 + m n + n^2)}{n^3} $$

Therefore, if $m\ne n$ we get that, $|\, a|<1$ and $\kappa_1=\sqrt{\frac{n(1 - |a|^2)}{m}}$.

In the case that $m=n$ a direct computation shows that  $p_3={\bf e_1}+\frac{2}{\kappa_1}\, {\bf e_2}+\frac{2}{\kappa_1}\, {\bf e_{m+3}}$ is a point in $M_f$, therefore $g(p_3)$ must be zero. It is not difficult to see that,

$$g(p_3)=\frac{64 (|\, a|^2-1 + \kappa_1^2)}{k_1} $$

Therefore, once again, when $m=n$, we get that $|\, a|<1$ and $\kappa_1=\sqrt{1-|\, a|^2}=\sqrt{\frac{n(1 - |\, a|^2)}{m}}$.

When $a\in \bfR{k}$ is the zero vector, the polynomial $f$ reduces to

$$  f(x)=2 x_1 x_2 +  \frac{m-n}{\sqrt{m n}} \, x_2^2  +\sqrt{\frac{n}{m}} \, |y|^2-\sqrt{ \frac{m}{n}} \, |z|^2  $$

When $a$ is not the zero vector,  doing a change of coordinates of the form,

$$\tilde{x}_i=x_i,\com{for $i=1,\dots m+n+2$ and} \tilde{x}_{m+n+3}=\frac{1}{|\, a|}\, (a_{n+m+1} x_{m+n+3}+\dots + a_{n+m+k}x_{n+m+k+2}) $$

we can assume, without loss of generality, that $a_{n+m+2}=\dots=a_{n+m+k+2}=0$. Under these assumptions the function $f$ reduces to,

$$f(x)=2 x_1 x_2 +  b  x_2^2 + \kappa_1 |y|^2 -\frac{m\kappa_1}{n}\, |y|^2+2 |\, a| x_2 x_{m+n+3} $$

Now, since $|\, a|<1$ we can find an isometric with respect to the metric $g_2$ that changes the coordinates so that

$$ \tilde{x}_1=\frac{x_1+|\, a|\, x_{m+n+3}}{\sqrt{1-|\, a|^2}},\quad  \tilde{x}_{m+n+3}=\frac{x_{m+n+3}+|\, a|\, x_1}{\sqrt{1-|\, a|^2}} \com{and} \tilde{x}_i=x_i\com{for all $i\ne1,m+n+3$}$$

After dividing $f$ by $\sqrt{1-|\, a|^2}$, the change of coordinates above will allow us to write  $f$ as

$$f(x)=2 x_1 x_2 + {b}  x_2^2 + {\kappa_1} |y|^2 -\frac{m{\kappa}_1}{n}\, |y|^2 $$

Once the polynomial $f$  has this form, by using the argument we used when $a$ is the zero vector, we can deduce that

$$  f(x)=2 x_1 x_2 +  \frac{m-n}{\sqrt{m n}} \, x_2^2  +\sqrt{\frac{n}{m}} \, |y|^2-\sqrt{ \frac{m}{n}} \, |z|^2  $$

This completes the proof of the theorem.
\end{proof}


\begin{thebibliography}{1}


\bibitem{C-Y} Cheng, S-Y., Yau, S-T \emph{Maximal space-like hypersurfaces in the Lorentz-Minkowski spaces}, Ann. Math. J. {\bf 104} 407-419 (1976)


\bibitem{H} Hsiang, W-Y \emph{Remarks on closed minimal submanifolds in the standard Riemannian m-sphere}, J. Differ. Geom. {\bf 1} 257-267 (1967)




\bibitem{Lawson} Lawson, B \emph{Complete minimal surfaces in $S^3$}, Ann. Math. {\bf 92} (2), (1970) 335-374.






\end{thebibliography}
\end{document}